\documentclass[review]{elsarticle}

\usepackage{hyperref}

\usepackage{amssymb}

\newtheorem{theorem}{Theorem}[section]

\journal{Mathematics and Computers in Simulation}

\usepackage{graphicx}
\usepackage{mathptmx}      
\usepackage{amssymb}

\usepackage{hyperref}
\usepackage{amsmath}
\usepackage{tikz}

\usepackage{soul}
\usepackage{xcolor}
\usepackage{amsfonts}

\usepackage[utf8]{inputenc}
\usepackage{bm}
\newcommand{\be}{\begin{equation}}
\newcommand{\ee}{\end{equation}}

\begin{document}

\begin{frontmatter}

\title{Asymptotically periodic and bifurcation points in fractional difference maps}

\author[mymainaddress,mysecondaryaddress]{Mark Edelman\corref{mycorrespondingauthor}}
\cortext[mycorrespondingauthor]{Corresponding author}
\ead{edelman@cims.nyu.edu}

\address[mymainaddress]{Stern College for Women, Yeshiva University, 245 Lexington Ave., New York, NY 10016}
\address[mysecondaryaddress]{Courant Institute of Mathematical Sciences at NYU, 251 Mercer Street, New York, NY 10012}

\begin{abstract}
The first step in investigating fractional difference maps, which do not
have periodic points except fixed points, is to find asymptotically periodic points and bifurcation points and draw asymptotic bifurcation diagrams. Recently derived equations that allow calculations of asymptotically periodic and bifurcation points contain coefficients defined as slowly converging infinite sums. In this paper we derive analytic expressions for coefficients of the equations that allow calculations of asymptotically periodic and bifurcation points in fractional difference maps.  
\end{abstract}


\end{frontmatter}


\section{Introduction}
\label{sec:1}

Fractional difference maps are maps with power-law-like memory. They are used to model biological (see, e.g. \cite{MEDie,Epi}) and socio-economic (see, e.g., \cite{TT1,TT2}) systems, memristors (see, e.g., \cite{Mem}), in image and signal encryption (see, e.g. \cite{I1,I2}), to control systems (see, e.g., \cite{C1,C2}), etc.

It is known that continuous and discrete fractional systems may
not have periodic solutions except fixed points (see, for example, \cite{PerD2,PerC2}). All bifurcation diagrams based on the finite time calculations on single trajectories are only approximations of
the asymptotic bifurcation diagrams. These approximations
depend on the initial conditions and the number of iterations.
But the asymptotically periodic solutions of fractional difference equations do exist, and the equations for finding these points in generalized fractional maps were derived in \cite{ME14,Helman,MELast}. These equations contain coefficients $S_{p,l}$ ($l$ is a period and $1 \le p \le l$) which are slowly converging series. The numerical evaluation of these series, in the case fractional and fractional difference maps, requires calculations of finite sums of tens of thousands of terms and calculations of the Riemann $\zeta$-function. It is also known, from the stability analysis of the discrete fractional systems (see \cite{MEstab}), that, in the case of fractional difference maps, the corresponding series may be summable (see, e.g., \cite{AbuSaris2013,Cermak2015,Mozyrska2015,Anh}). The equations that define bifurcation points of fractional difference maps \cite{MEBif} also depend on the same coefficients (sums). 

Following this introduction, in the preliminaries (Section \ref{sec:2}), based on well-known results, we show that all aspects of the asymptotic theory of generalized fractional maps, such as periodic points, bifurcation points, and conditions of stability, strongly depend on the values of $S_{p,l}$ (see also a recent review \cite{MEBor}). Then, in Section \ref{sec:3}, using a recently derived combinatorial identity (see \cite{MSS}), we derive analytic expressions for the coefficients (sums) of the equations defining asymptotically periodic and bifurcation points in the case of fractional difference maps. In Section \ref{sec:4}, we discuss the validity of our results, their fundamental significance, computational aspects of calculations of asymptotically periodic and bifurcation points of fractional difference maps, and possible applications. The concluding remarks are presented in Section \ref{sec:5}. The Appendix contains C-codes for calculations of sums $S_{p,l}$ based the analytic formulae derived in this paper and on the numerical algorithm proposed in \cite{ME14} (Eq.~(35) in that paper).

\section{Preliminaries}
\label{sec:2}

For $0<\alpha<1$, the generalized universal $\alpha$-family of maps 
is defined as (see \cite{ME14,Helman}):
\begin{eqnarray}
x_{n}= x_0 
-\sum^{n-1}_{k=0} G^0(x_k) U_\alpha(n-k),
\label{FrUUMapN}
\end{eqnarray} 
where $G^0(x)=h^\alpha G_K(x)/\Gamma(\alpha)$, $x_0$ is the initial condition, $h$ is the time step of the map, $\alpha$ is the order of the map, $G_K(x)$ is a nonlinear function depending on the parameter $K$, $U_\alpha(n)=0$ for $n \le 0$, and $U_\alpha(n) \in \mathbb{D}^0(\mathbb{N}_1)$. The space $\mathbb{D}^i(\mathbb{N}_1)$ is defined as (see \cite{Helman})
{\setlength\arraycolsep{0.5pt}
\begin{eqnarray}
&&\mathbb{D}^i(\mathbb{N}_1)\ \ = \ \ \{f: |\sum^{\infty}_{k=1}\Delta^if(k)|>N, \ \ \forall N, \ \ N \in \mathbb{N}, 
\sum^{\infty}_{k=1}|\Delta^{i+1}f(k)|=C, \ \ C \in \mathbb{R}_+\},
\label{DefForm}
\end{eqnarray}
}
where $\Delta$ is a forward difference operator defined as
\begin{equation}
\Delta f(n)= f(n+1)-f(n).
\label{Delta}
\end{equation}
In the case Caputo fractional difference maps, which are defined as solutions of the Caputo h-difference equation \cite{DifSum,Fall,Chaos2014}
\begin{equation}
(_0\Delta^{\alpha}_{h,*} x)(t) = -G_K(x(t+(\alpha-1)h)),
\label{LemmaDif_n_h}
\end{equation}
where $t\in (h\mathbb{N})_{m}$, with the initial conditions 
 \begin{equation}
(_0\Delta^{k}_h x)(0) = c_k, \ \ \ k=0, 1, ..., m-1, \ \ \ 
m=\lceil \alpha \rceil,
\label{LemmaDifICn_h}
\end{equation}
the kernel $U_\alpha(n)$ is the falling factorial function: 
{\setlength\arraycolsep{0.5pt}
\begin{eqnarray}
&&U_{\alpha}(n)=(n+\alpha-2)^{(\alpha-1)} 
, \   \ U_{\alpha}(1)=(\alpha-1)^{(\alpha-1)}=\Gamma(\alpha).
\label{UnFrDif}
\end{eqnarray} 
}
The definition of the falling factorial $t^{(\alpha)}$ is
\begin{equation}
t^{(\alpha)} =\frac{\Gamma(t+1)}{\Gamma(t+1-\alpha)}, \ \ t\ne -1, -2, -3.
...
\label{FrFacN}
\end{equation}
The falling factorial is asymptotically a power function:
\begin{equation}
\lim_{t \rightarrow
  \infty}\frac{\Gamma(t+1)}{\Gamma(t+1-\alpha)t^{\alpha}}=1,  
\ \ \ \alpha \in  \mathbb{R}.
\label{GammaLimitN}
\end{equation}
The $h$-falling factorial $t^{(\alpha)}_h$ is defined as
\begin{eqnarray}
&&t^{(\alpha)}_h =h^{\alpha}\frac{\Gamma(\frac{t}{h}+1)}{\Gamma(\frac{t}{h}+1-\alpha)}= h^{\alpha}\Bigl(\frac{t}{h}\Bigr)^{(\alpha)}, \  \
\frac{t}{h} \ne -1, -2, -3,
....
\label{hFrFacN}
\end{eqnarray}
Majority of the introduced, investigated, and used in applications maps are Caputo fractional difference maps.

The following equations define period-$l$ points in generalized fractional maps of the orders $0<\alpha<1$ \cite{ME14}:
{\setlength\arraycolsep{0.5pt}   
\begin{eqnarray} 
&&x_{lim,m+1}-x_{lim,m}=S_{1,l} G^0(x_{lim,m})+\sum^{m-1}_{j=1}S_{j+1,l} G^0(x_{lim,m-j})
\nonumber \\
&&+\sum^{l-1}_{j=m}S_{j+1,l} G^0(x_{lim,m-j+l}), \  \ 0<m<l,
\label{LimDifferences}
\\
&&\sum^{l}_{j=1} G^0(x_{lim,j})=0,
\label{LimDifferencesN}
\end{eqnarray}
}
where
{\setlength\arraycolsep{0.5pt}   
\begin{eqnarray} 
&&S_{j+1,l}=\sum^{\infty}_{k=0}\Bigl[
U_{\alpha} (lk+j) - U_{\alpha} (lk+j+1)\Bigr], \  \ 0 \le j<l.
\label{Ser}
\end{eqnarray}
}
It is easy to see that 
\begin{equation}
\sum^{l}_{j=1}S_{j,l}=0.
\label{Ssum}
\end{equation}

In the case of $p$-dimensional maps ($1 \le i \le p$) (see \cite{MELast})
\begin{eqnarray}
x_{i,n} = x_{i,0}  
-\sum^{n-1}_{k=0} {G}_i^0(x_{1,k},x_{2,k},..., x_{p,k}) U_{\alpha_i}(n-k),
\label{FrUUMapNM}
\end{eqnarray} 
the periodic points are defined as solutions of the system of $(l-1)\times p$ equations:
{\setlength\arraycolsep{0.5pt}   
\begin{eqnarray} 
&&x_{i,l,m+1}-x_{i,l,m}= \sum^{m-1}_{j=0}S_{i,j+1,l} G_i^0(x_{1,l,m-j},x_{2,l,m-j}, ..., x_{p,l,m-j})
\nonumber \\
&&+\sum^{l-1}_{j=m}S_{i,j+1,l} G_i^0(x_{1,l,m-j+l},x_{2,l,m-j+l}, ..., x_{p,l,m-j+l}), 
\nonumber \\
&&0<m<l, \  \ 0<i \le p
\label{LimDifferencesNM}
\end{eqnarray}
}
and additional $p$ equations
\begin{equation}
\sum^{l}_{j=1} G_i^0(x_{1,l,j},x_{2,l,j}, ..., x_{p,l,j})=0, \ 0 < i \le p.
\label{Close}
\end{equation}

Bifurcation points in the maps of the order $0<\alpha<1$ are defined by the Theorem 1 from \cite{MEBif}:
\begin{theorem}\label{Th1} 
The $T=2^{n-1}$ --  $T=2^{n}$ bifurcation points, $2^{n-1}$ values of $x_{2^{n-1}bif,i}$ with $0<i \le 2^{n-1}$ and the value of the nonlinear parameter $K_{2^{n-1}bif}$, of a fractional generalization of a nonlinear one-dimensional map $x_{n+1}=F_K(x_n)$ written as the Volterra difference equations of convolution type
\begin{eqnarray}
x_{n}= x_0 
-\sum^{n-1}_{k=0} G^0(x_k) U_\alpha(n-k),
\label{FrUUMapN1}
\end{eqnarray} 
where $G^0(x)=h^\alpha G_K(x)/\Gamma(\alpha)$, $x_0$ is the initial condition, $h$ is the time step of the map, $\alpha$ is the order of the map, $G_K(x)=x-F_K(x)$, $U_\alpha(n)=0$ for $n \le 0$, $U_\alpha(n) \in \mathbb{D}^0(\mathbb{N}_1)$, and
\begin{eqnarray}
&&\mathbb{D}^i(\mathbb{N}_1) \ \ = \ \ \left\{f: \left|\sum^{\infty}_{k=1}\Delta^if(k)\right|>N, 
\right.\nonumber \\
&& \left.\forall N, \ \ N \in \mathbb{N}, 
\sum^{\infty}_{k=1}\left|\Delta^{i+1}f(k)\right|=C, \ \ C \in \mathbb{R}_+\right\},
\label{DefForm1}
\end{eqnarray}
are defined by the system of $2^{n-1}+1$ equations   
\begin{eqnarray} 
&&x_{2^{n-1}bif,m+1}-x_{2^{n-1}bif,m}=S_{1, 2^{n-1}}    
G^0(x_{2^{n-1}bif,m})\nonumber\\
&&\hspace{25pt}{}+\sum^{m-1}_{j=1}S_{j+1,2^{n-1}} G^0(x_{2^{n-1}bif,m-j})\nonumber \\
&&\hspace{25pt}{}+\sum^{2^{n-1}-1}_{j=m}S_{j+1,2^{n-1}} G^0(x_{2^{n-1}bif,m-j+2^{n-1}}),\label{LimDifferencesNN}\\
&&\hspace{150pt} 0<m<2^{n-1},\nonumber
\\
&&\sum^{2^{n-1}}_{j=1} G^0(x_{2^{n-1}bif,j})=0,
\label{LimDifferencesNN2}
\\
&& \operatorname{det}(A)=0
\label{DetNN},
\end{eqnarray}
where  
\begin{eqnarray} 
&&S_{j+1,l}=\sum^{\infty}_{k=0}\Bigl[
U_{\alpha} (lk+j) - U_{\alpha} (lk+j+1)\Bigr], \nonumber\\
&&\hspace{50pt} 0 \le j<l, \quad
S_{i,l}=S_{i+l,l}, \quad i \in \mathbb{Z},
\label{SerNN1}
\end{eqnarray}
and the elements of the $2^{n-1}$-dimensional matrix $A$ are 
\begin{eqnarray} 
&&\hspace{-10pt}A_{i,j}=\frac{dG^0(x)}{dx}\Bigl|_{x_{2^{n-1}bif,j}}
 \sum^{i+2^{n-1}-1}_{m=i}S_{m-j+1,2^n}+\delta_{i,j}, 
 \label{DetNN2} 
 \\
 &&\hspace{140pt} 0 < i,j \le 2^{n-1}. \nonumber
\end{eqnarray}

\end{theorem}

\section{Sums $S_{p,l}$ for $l$-cycles of fractional difference maps}
\label{sec:3}

As one may see from Section~\ref{sec:2}, all aspects of the asymptotic theory of fractional difference maps depend on the values of the slowly converging sums $S_{p,l}$. As it is shown in \cite{MEstab}, the stability conditions of the map's fixed points also depend on $S_{p,2}$. This is why obtaining analytic expressions for $S_{p,l}$ is an important part of the asymptotic theory of the fractional difference maps.

The definition of $S_{p,l}$ from \cite{ME14} for fractional difference maps may be rewritten using the following chain of transformations:
{\setlength\arraycolsep{0.5pt}   
\begin{eqnarray} 
&&S_{p,l}= 
\sum^{\infty}_{k=0}\Bigl[(lk+p+\alpha-3)^{(\alpha-1)}- (lk+p+\alpha-2)^{(\alpha-1)} \Bigr] =\sum^{\infty}_{k=0}
\Bigl[\frac{\Gamma(lk+p+\alpha-2)}{ \Gamma(lk+p-1)}
\nonumber \\
&&- \frac{\Gamma(lk+p+\alpha-1)}
{ \Gamma(lk+p)} \Bigr] =(1-\alpha)\sum^{\infty}_{k=0}
\frac{\Gamma(lk+p+\alpha-2)}{\Gamma(lk+p)}
=-\Gamma(\alpha) \sum^{\infty}_{k=0}
\frac{\Gamma(lk+p+\alpha-2)}{\Gamma(\alpha-1)\Gamma(lk+p)}
\nonumber \\
&&=-\Gamma(\alpha) \sum^{\infty}_{k=0}
\left( \begin{array}{c}
lk+p+\alpha-3 \\ lk+p-1
\end{array} \right)
=\Gamma(\alpha) (-1)^{p}\sum^{\infty}_{k=0}(-1)^{lk}
\left( \begin{array}{c}
1-\alpha \\ lk+p-1
\end{array} \right).
\label{Sp_l}
\end{eqnarray}
}
Using absolute convergence of series and the following identity (see \cite{MSS})
{\setlength\arraycolsep{0.5pt}   
\begin{eqnarray} 
&&\sum^{\infty}_{k=0}
\left( 
\begin{array}{c}
\gamma \\ t+ks
\end{array} \right)
= \frac{1}{s}\sum^{s-1}_{j=0} \omega^{-jt}(1+\omega^j)^\gamma,
\label{I1_3}
\end{eqnarray}
}
where $\omega=e^{i2\pi/s}$, for the even and odd periods we obtain
{\setlength\arraycolsep{0.5pt}   
\begin{eqnarray} 
&&S_{p,2n}= \Gamma(\alpha) (-1)^{p}\sum^{\infty}_{j=0}
\left( \begin{array}{c}
1-\alpha \\ 2nj+p-1
\end{array} \right)
=\frac{\Gamma(\alpha)}{2n}(-1)^p\sum^{2n-1}_{j=0}e^{-i\pi j(p-1)/n}(1+ e^{i\pi j/n})^{1-\alpha}
\nonumber \\
&&=\frac{\Gamma(\alpha)}{2n}(-1)^p\Bigl[2^{1-\alpha}
+\sum^{n-1}_{j=1}\Bigl(e^{-i\pi j(p-1)/n} e^{i\pi j(1-\alpha)/(2n)}(2\cos(\pi j/(2n)))^{1-\alpha}
\nonumber \\
&&+ e^{i\pi j(p-1)/n} e^{-i\pi j(1-\alpha)/(2n)} (2\cos(\pi j/(2n)))^{1-\alpha}\Bigr) \Bigr]
\nonumber \\
&&=\frac{\Gamma(\alpha) 2^{-\alpha}}{n}(-1)^p\Bigl[1
+2\sum^{n-1}_{j=1}(\cos(\pi j/(2n)))^{1-\alpha}\cos(\pi j(2p+\alpha-3)/(2n)) \Bigr],
\label{Sp_2n}
\end{eqnarray}
}
{\setlength\arraycolsep{0.5pt}   
\begin{eqnarray} 
&&S_{p,2n+1}= \Gamma(\alpha) (-1)^{p}\sum^{\infty}_{j=0}(-1)^j
\left( \begin{array}{c}
1-\alpha \\ (2n+1)j+p-1
\end{array} \right)
\nonumber \\
&&=\Gamma(\alpha) (-1)^{p} \Bigl\{ \sum^{\infty}_{j=0}
\left( \begin{array}{c}
1-\alpha \\ 2(2n+1)j+p-1
\end{array} \right)
-\sum^{\infty}_{j=0}
\left( \begin{array}{c}
1-\alpha \\ 2(2n+1)j+p+2n
\end{array} \right)
\Bigr\}
\nonumber \\
&&=\frac{\Gamma(\alpha)}{2(2n+1)}(-1)^p \Bigl\{   \sum^{4n+1}_{j=0}e^{-\frac{i\pi j(p-1)}{2n+1}}(1+ e^{\frac{i\pi j}{2n+1}})^{1-\alpha}-\sum^{4n+1}_{j=0}e^{-\frac{i\pi j(p+2n)}{2n+1}}(1
\nonumber \\
&&+ e^{\frac{i\pi j}{2n+1}})^{1-\alpha}\Bigr\}
=\frac{\Gamma(\alpha)}{2(2n+1)}(-1)^p    \sum^{2n}_{j=1} \Bigl\{      \Bigl(e^{-\frac{i\pi j(p-1)}{2n+1}}- e^{-\frac{i\pi j(p+2n)}{2n+1}}\Bigr)     (1
\nonumber \\
&&+ e^{\frac{i\pi j}{2n+1}})^{1-\alpha}
+\Bigl(e^{\frac{i\pi j(p-1)}{2n+1}}- e^{\frac{i\pi j(p+2n)}{2n+1}}\Bigr)     (1+ e^{-\frac{i\pi j}{2n+1}})^{1-\alpha}\Bigr\}
\nonumber \\
&&=\frac{\Gamma(\alpha)}{2n+1}(-1)^p i    \sum^{2n}_{j=1}\Bigl(2\cos\frac{\pi j}{2(2n+1)} \Bigr)^{1-\alpha}
\sin\frac{j\pi}{2}\times
\Bigl(e^{-\frac{i\pi j(2p+2n-2+\alpha)}{2(2n+1)}}
\nonumber \\
&&- e^{\frac{i\pi j(2p+2n-2+\alpha)}{2(2n+1)}}\Bigr)
=\frac{2^{2-\alpha}\Gamma(\alpha)}{2n+1}(-1)^p   
\sum^{2n}_{j=1}\Bigl(\cos\frac{\pi j}{2(2n+1)} \Bigr)^{1-\alpha}
\sin\frac{j\pi}{2}
\nonumber \\
&&\times
\sin\frac{\pi j(2p+2n-2+\alpha)}{2(2n+1)}
\nonumber \\
&&=\frac{2^{2-\alpha}\Gamma(\alpha)}{2n+1}(-1)^p   
\sum^{n-1}_{j=0}\Bigl(\cos\frac{\pi (2j+1)}{2(2n+1)} \Bigr)^{1-\alpha}
(-1)^j
\nonumber \\
&&\times
\sin\frac{\pi (2j+1)(2p+2n-2+\alpha)}{2(2n+1)}
\label{Sp_2n+1}
\end{eqnarray}
}

\section{Numerical simulations and discussion}
\label{sec:4}

Although in \cite{MSS}, Eq.~(\ref{I1_3}) is proven for integer values of $\gamma \in \mathbb{N}_0$, it is valid for any real values $0 \le \gamma \le 1$. To verify this, we calculated the values of $S_{i,4}$ ($1 \le i \le 4$) using Eq.~(\ref{Sp_2n}) for $\alpha=0.5$ and $\alpha=0.99$. A C-code used in our calculations is the second code presented in the Appendix.

This is the output of the C-code utilizing Eq.~(\ref{Sp_2n})  for $\alpha=0.5$:

-1.600340512580117, 1.029969956588537, 0.347026375264616, 0.223344180726962, 

5.000000e-01, -0.000000000000001  \\
and for $\alpha=0.99$:

-0.757883584649114, 0.257180815994029, 0.251449398374527, 0.249253370280558, 
   
9.900000e-01, -0.000000000000000.

The first four values of the output are $S_{1,4}$, $S_{2,4}$, $S_{3,4}$, and $S_{4,4}$, the fifth value is the value of the order $\alpha$, and the last one is the total $S_{1,4}+S_{2,4}+S_{3,4}+S_{4,4}$, which theoretical value is equal to zero.
The results are compared to the corresponding values obtained using the expression which allows a fast calculation of the series using tens of thousands of operations (see Eq.~(35) in \cite{ME14}).  The corresponding C-code is the first code presented in the Appendix.

The output of the code based on the summation algorithm Eq.~(35) from \cite{ME14} for $\alpha=0.5$:

-1.600340512580117, 1.029969956588546, 0.347026375264620, 0.223344180726963, 
 
5.000000e-01, 0.000000000000012  \\
and for $\alpha=0.99$:

-0.757883584649062, 0.257180815994081, 0.251449398374580, 0.249253370280609, 
   
9.900000e-01, 0.000000000000208.

The computational time on a single processor of the one of the Courant's computers is less than $10^{-3}$s for both codes. The relative accuracy of the calculations in the former case is $10^{-16}$ and in the latter case it is $10^{-14}$. The lower accuracy of the calculations using the summation algorithm apparently is a result of hundreds of thousands of operations required to complete calculations. It is obvious that both algorithms provide sufficient speed and accuracy for most practical cases.

The computation of the coefficients $S_{i,l}$ is a part of the algorithm of the computation of periodic and bifurcation points used in \cite{ME14,MEBif}, where it was applied to investigate the fractional logistic maps (fractional extensions of the map $x_{n+1}=Kx_n(1-x_n)$). It took less than a second to calculate numerically all $S_{i,l}$ values. Computations of the periodic points were significantly more time-consuming. They required solving systems of the thousands of second order equations using standard Mathematica and Matlab algorithms. These algorithms converge only when the initial approximations to the solutions are very close to the real solutions. To draw the bifurcation diagrams, the calculations were performed with a small step in the nonlinear parameter $K$ using the solutions from the previous step as the initial conditions for the next step. An additional time-consuming task was the separation of the stable periodic points from unstable lower-period points. It took a couple of months to calculate the data for the bifurcation diagram Fig.~1 and Table~1 from \cite{MEBif}.
So, one may see that the analytic formulae for $S_{i,l}$ are important but they do not play a critical role in drawing the bifurcation diagrams. 
 
The analytic expressions for $S_{i,l}$ obtained in this paper may be used to obtain analytic expressions for low-period periodic and bifurcation points where it is possible. This may be important for the analysis of the low-period behavior of various discrete systems. 

Another area of application of the results obtained in this paper is the analysis of the universality in fractional dynamics (for the universality in regular dynamics see, e.g., \cite{Cvi,Fei}). The results of the computations presented in Tables~1~and~2 from \cite{MEBif} show that the ratios of the parameter intervals between the consecutive asymptotic bifurcations in fractional maps converge to the Feigenbaum's constant $\delta$ \cite{FC} but much slower than in the integer case. The proof of universality in fractional case is more complicated for the following reasons: a) the periodicity in fractional case exists only in the asymptotic sense, and the universality should be only the asymptotic universality; b) fractional maps are maps with memory, and their Poincar$\acute{e}$ plots (graphs $x_{n+1}$ vs. $x_{n}$) used by Feigenbaum in his analysis of universality have a limited value. Obtaining analytic expressions for the coefficients of the equations defining asymptotic bifurcation points is a small but important step in the proof of the asymptotic universality in fractional difference maps (discrete systems with the falling factorial memory).

\section{Conclusion}
\label{sec:5}

In this paper we derived the analytic expressions for the coefficients of the equations that define periodic and bifurcation points in fractional difference maps of the orders $0<\alpha<1$ (Eqs.~(\ref{Sp_2n})~and~(\ref{Sp_2n+1})). 
The purpose of this publication was not to investigate any particular fractional difference map, but rather to provide the framework which will enable researchers to investigate any maps with fractional factorial memory. Our results will allow researchers to draw asymptotic bifurcation diagrams of fractional difference maps by solving Eqs.~(\ref{LimDifferences})~and~(\ref{LimDifferencesN}) in the case of maps of the orders $0<\alpha<1$, or Eqs.~(\ref{LimDifferencesNM})~and~(\ref{Close}) in the case of multidimensional fractional maps. Calculations of coefficients (sums) $S_{p,l}$ of these equations should be a part of the corresponding numerical algorithms.  

Using analytic expressions derived in this paper instead of adding tens of thousands of terms based on Eq.~(35) from \cite{ME14} will make calculations of the coefficients thousands of times faster, more accurate, and will make it unnecessary to use Tables~2~and~4 from \cite{ME14} for finding $S_{p,l}$.

We should note two outstanding problems related to this paper to be addressed in further publications: 

1. In this paper, the validity of the identity Eq.~(\ref{I1_3}) for fractional values of $\gamma$ is demonstrated numerically, but the theoretical proof is still due.

2. Based on the results of numerical simulations, in \cite{MEBif}, the authors made a conjecture that the Feigenbaum number $\delta$ exists in fractional difference maps and has the same value as in regular maps. Derivation of the analytic expression for $S_{p,l}$ in this paper is a small step which, in conjunction with Theorem~1 from \cite{MEBif} defining bifurcation points, may lead to a theoretical proof of this conjecture.

\section*{Appendix}
\label{sec:A}

C-codes used to calculate $S_{p,4}$ based on Eq.~(35) from \cite{ME14}: \newline
\#include $<$stdio.h$>$ \\
\#include $<$math.h$>$  \\
\#include $<$stdlib.h$>$ \\
\#define NIter 19999.1 \\  
\#define l 4 \\
double Zeta(double s); \\
double gamma(double xx); \\
main(int argc, char $*$argv[])  \\
$\{$   \\                                  
\hspace*{.2cm} double  x, sum, S[l], SZ[l], St, S1, S2, S3, alp, beta, Z2, Z3, Z4, g1, Ssum; \\
\hspace*{.2cm} int j, k, m; \\
\hspace*{.2cm}  FILE  $*$in, $*$out;        \\
\hspace*{.2cm}   if(argc $>$ 2) $\{$    \\
 \hspace*{.4cm}   if( (in = fopen(argv[1],"r")) == NULL) $\{$ \\
 \hspace*{.6cm}    printf("Can't open file \%s\textbackslash n", argv[1]);  \\
\hspace*{.6cm}      exit(1);      $\}$            \newline
\hspace*{.4cm}       fscanf(in, "\%le", \&alp);         \newline    
\hspace*{.4cm}       fclose(in);    $\}$                    \newline
\hspace*{.2cm}     if( (out = fopen(argv[2],"w")) == NULL) $\{$   \newline
\hspace*{.4cm}    printf("Can't open file \%s\textbackslash n", argv[2]);         \\
\hspace*{.4cm}    exit(1);     $\}$    \\ 
\hspace*{.2cm}  k=0; \\
\hspace*{.2cm}  Ssum=0; \\
\hspace*{.2cm}   S1=S2=S3=0.; \\
\hspace*{.2cm}   while(k$<$NIter) \\
\hspace*{.4cm}    $\{$ \\
\hspace*{.4cm}      k++; \\
\hspace*{.4cm}      S1+=exp((alp-2.)*log(k)); \\
\hspace*{.4cm}     S2+=exp((alp-3.)*log(k));  \\
\hspace*{.4cm}      S3+=exp((alp-4.)*log(k)); \\
\hspace*{.4cm}     $\}$  \\
\hspace*{.2cm}  Z2=Zeta(2-alp)-S1; \\
\hspace*{.2cm}  Z3=Zeta(3-alp)-S2; \\
\hspace*{.2cm}  Z4=Zeta(4-alp)-S3; \\
\hspace*{.2cm} beta=1.-alp;  \\
\hspace*{.2cm}  g1=gamma(alp-1.); \\
\hspace*{.2cm}   for(j=0; j$<$l; j++)  \\
\hspace*{.4cm}    $\{$ \\
\hspace*{.4cm}    St=g1; \\
\hspace*{.4cm}     for(k=0; k$<$j; k++) \\
\hspace*{.6cm}      St=St*(alp+k-1)/(k+1); \\
\hspace*{.4cm}    S[j]=St; \\
\hspace*{.4cm}     k=0.; \\
\hspace*{.4cm}     while(k$<$NIter) \\
\hspace*{.6cm}      $\{$ \\
\hspace*{.6cm}      k++;  \\
\hspace*{.6cm}       for(m=0; m$<$l; m++) \\
\hspace*{.8cm}          St=St*(l*(k-1)+j+alp-1+m)/(l*(k-1)+j+1+m); \\
\hspace*{.6cm}        S[j]+=St; \\
\hspace*{.6cm}        $\}$ \\
\hspace*{.4cm}    S[j]=beta*(S[j]+ exp((alp-2.)*log(l))*(Z2
    +0.5*(alp-2.)/l*((2*j+alp-1.)*Z3+ \\
\hspace*{.4cm}     (alp-3.)*(3.*(2*j+alp-1)*(2*j+alp-1)-alp+1.)*Z4/(12.*l)))); \\
\hspace*{.4cm}     fprintf(out, "\%17.15f \textbackslash n", S[j]); \\
\hspace*{.4cm}   Ssum=Ssum+S[j];  \\
\hspace*{.4cm}   $\}$  \\
\hspace*{.2cm}      fprintf(out, "\%le, \%17.15f  \textbackslash n", alp, Ssum); \\
\hspace*{.2cm}  fclose(out); \\
\hspace*{.2cm} $\}$

and on Eq.~(\ref{Sp_2n}) derived in this paper: \\
\#include $<$stdio.h$>$ \newline
\#include $<$math.h$>$   \newline
\#include $<$stdlib.h$>$ \newline
\#define n 2 \newline
double gamma(double xx);  \newline
main(int argc, char $*$argv[])  \newline
$\{$                                     \newline
\hspace*{.2cm}  double x, sum, S[2*n], alp, g1, Ssum;  \newline
\hspace*{.2cm}  int i, j, m;             \newline
\hspace*{.2cm}  FILE  $*$in, $*$out;        \newline
\hspace*{.2cm}   if(argc $>$ 2) $\{$    \newline
 \hspace*{.4cm}   if( (in = fopen(argv[1],"r")) == NULL) $\{$ \newline
 \hspace*{.6cm}    printf("Can't open file \%s\textbackslash n", argv[1]);  \newline
\hspace*{.6cm}      exit(1);      $\}$            \newline
\hspace*{.4cm}       fscanf(in, "\%le", \&alp);         \newline    
\hspace*{.4cm}       fclose(in);    $\}$                    \newline
\hspace*{.2cm}     if( (out = fopen(argv[2],"w")) == NULL) $\{$   \newline
\hspace*{.4cm}    printf("Can't open file \%s\textbackslash n", argv[2]);         \newline 
\hspace*{.4cm}    exit(1);     $\}$                            \newline
\hspace*{.2cm}  g1=gamma(alp)* exp((-alp)*log(2.))/n;  \\
\hspace*{.2cm}   j=1; \\
\hspace*{.2cm}   Ssum=0; \\
\hspace*{.2cm}   for(i=0; i$<$2*n; i++) \\
\hspace*{.4cm}   $\{$  \\
\hspace*{.4cm}    j=-j;  \\
\hspace*{.4cm}    S[i]=1;  \\
\hspace*{.4cm}    for(m=1; m$<$n; m++)  \\
\hspace*{.6cm}      $\{$   \\
\hspace*{.6cm}       S[i]=S[i]+2.0*exp((1.-alp)*log(cos(3.141592653589793*m/2.0/n)))  \\
 \hspace*{.6cm}  *cos(3.141592653589793*m*(2.0*(i+1)+alp-3)/2.0/n);  \\
\hspace*{.6cm}   $\}$   \\
\hspace*{.4cm}      S[i]=g1*j*S[i];  \\
\hspace*{.4cm}     fprintf(out, "\%17.15f \textbackslash n", S[i]); \\
\hspace*{.4cm}      $\}$ \\
\hspace*{.4cm}   Ssum=Ssum+S[i];  \\
\hspace*{.2cm}      fprintf(out, "\%le, \%17.15f  \textbackslash n", alp, Ssum); \\ 
\hspace*{.2cm}      fclose(out);  \\
    $\}$

Here Zeta(double s) and gamma(double xx) are some standard algorithms for calculations of the Riemann $\zeta$-function and the Gamma function freely available on the Internet. They are also available on request from the author.

\section*{Data availability}

No data was used for the research described in the article.

\section*{Acknowledgements}
The author acknowledges support from Yeshiva University and expresses his gratitude to the administration of Courant Institute of Mathematical Sciences at NYU for the opportunity to complete all computations at
Courant.



\begin{thebibliography}{10}


\bibitem{MEDie} M. Edelman, Evolution of systems with power-law memory: Do we have to die? (Dedicated to the Memory of Valentin Afraimovich), in: C.H. Skiadas, C. Skiadas (Eds.), Demography of Population Health, Aging and Health Expenditures, pp. 65--85. Springer, eBook, 2020.


 
\bibitem{Epi} A.G.M. Selvam, D.A. Vianny, Discrete fractional order SIR epidemic model of childhood diseases with constant vaccination and its stability. Int. J. Tech. Innova. Mod. Engin. Sci. 4 (2018) 405--410. 

\bibitem{TT1}
V.E. Tarasov, V.V. Tarasova, Economic Dynamics with Memory: Fractional Calculus Approach, De Gruyter, Berlin, Boston, 2021. 

\bibitem{TT2} Tarasov, V.E., V.V. Tarasova, Logistic map with memory from economic model. Chaos. Soliton. Fract. 95 (2017) 84--91. 

\bibitem{Mem} S. He, K. Sun, Y. Peng, L. Wang, Modeling of discrete fracmemristor and its application. AIP Adv. 10 (2020) 015332.

\bibitem{I1}L.-L. Huang, D. Baleanu, G.-C. Wu, S.-D. Zeng, A new application of the fractional logistic map. Rom. J. Phys. 61, (2016) 1172--1179. 

\bibitem{I2}
D. Ding, J. Wang, M. Wang, et al., Controllable multistability of fractional-order memristive coupled chaotic map and its application in medical image encryption. Eur. Phys. J. Plus 138 (2023) 908. 


\bibitem{C1} P. Ostalczyk, Discrete Fractional Calculus: Applications in Control and Image Processing, World Scientific, River Edge, 2016.


\bibitem{C2} M. Ortigueira, Discrete-time fractional difference calculus: origins, evolutions, and new formalisms, Fractal Fract. 7 (2023) 502. 






\bibitem{PerD2} J.~Jagan~Mohan, Quasi-periodic solutions of fractional nabla difference systems, 
Fract. Differ. Calc. 7 (2017) 339--355. 


\bibitem{PerC2} E.~Kaslik, S.~Sivasundaram,
Nonexistence of periodic solutions in fractional order
dynamical systems and a remarkable difference bet
ween integer and fractional order derivatives of periodic functions,
Nonlinear Anal. Real World Appl. 13 (2012) 1489--1497. 

\bibitem{ME14} M.~Edelman,
Cycles in asymptotically stable and chaotic fractional maps,
Nonlin. Dyn. 104 (2021) 2829--2841. 


\bibitem{Helman}
M.~Edelman, A.B.~Helman, Asymptotic cycles in fractional maps of arbitrary positive orders, Fract. Calc. Appl. Anal. (2022) https://doi.org/10.1007/s13540-021-00008-w

\bibitem{MELast} M.~Edelman,
Asymptotic cycles in fractional generalizations of
Multidimensional maps, Fract. Calc. Appl. Anal. (2025)
https://doi.org/10.1007/s13540-024-00364-3

\bibitem{MEstab}
M.~Edelman,
Stability of Fixed Points in Generalized Fractional Maps of the Orders $0< \alpha <1$, Nonlin. Dyn. 111 (2023) 10247--10254. https://doi.org/10.1007/s11071-023-08359-0.


\bibitem{AbuSaris2013}
R.~Abu-Saris, Q.~Al-Mdallal, On the asymptotic stability
of linear system of fractional-order difference equations, Fract.
Calc. Appl. Anal. 16 (2013) 613--629.

\bibitem{Cermak2015}
J.~\v{C}erm\'{a}k, I.~Gy\H{o}ri, L.~Nechv\'{a}tal, On explicit stability conditions for a linear fractional difference system, Fract. Calc. and Appl. Anal. 18 (2015) 651--672.

\bibitem{Mozyrska2015}
D.~Mozyrska, M.~Wyrwas, The z-transform method and delta type fractional difference operators, Discrete Dyn. Nat. Soc. 2015 (2015) 852734.


\bibitem{Anh} P.T.~Anh, A.~Babiarz, A.~Czornik, M.~Niezabitowski,  S.~Siegmund, Asymptotic properties of discrete linear fractional equations, B. Pol. Acad. Sci.-Tech. 67 (2019) 749--759.

\bibitem{MEBif} M.~Edelman, A.B.~Helman, R.~Smidtaite,
Bifurcations and transition to chaos in generalized fractional
maps of the orders $0<\alpha<1$, Chaos 33 (2023) 063123. https://doi.org/10.1063/5.0151812

\bibitem{MEBor} M.~Edelman,
Periodic Points, Stability, Bifurcations, and Transition to Chaos in Generalized Fractional Maps,
IFAC-PapersOnLine, Vol. 58, Is. 12 (2024) Pages 131--142.

\bibitem{MSS} A.T.~Benjamin, B.~Chen, K.~Kindred, Sums of evenly spaced binomial coefficients,
Maths. Mag. 83 (2010) 370--373.


\bibitem{DifSum} F.~Chen, X.~Luo, Y.~Zhou,
Existence Results for Nonlinear Fractional Difference Equation,
Adv. Differ. Eq. 2011 (2011) 713201.

\bibitem{Fall} G.-C.~Wu, D.~Baleanu, S.-D.~Zeng,
Discrete chaos in fractional sine and standard maps,
Phys. Lett. A 378 (2014) 484--487.

\bibitem {Chaos2014} M.~Edelman,
Caputo standard $\alpha$-family of maps: 
Fractional difference vs. fractional,
Chaos 24 (2014) 023137.

\bibitem{Cvi} P. Cvitanovic, Universality in Chaos, Taylor \& Francis Group, New York, 1989.


\bibitem{Fei} M.J. Feigenbaum, Quantitative universality for a class of nonlinear transformations, J. Stat. Phys. 19 (1978) 25--52.


\bibitem{FC}
Feigenbaum constants, Wikipedia
https://en.wikipedia.org/wiki/Feigenbaum\_constants.
















\end{thebibliography}
\end{document}